\documentclass{article}
\usepackage{acra}

\usepackage[T1]{fontenc}
\usepackage[utf8]{inputenc}
\usepackage{enumerate,amsmath,amsthm}
\usepackage{amssymb}
\usepackage{bookman}
\usepackage{enumitem,etoolbox,hyperref}
\usepackage{eulervm} % Neuvěřitelný font!!!!!
\usepackage{tabu}
\usepackage{multirow}
\usepackage{graphicx} %% Obrázky
\usepackage{wrapfig} %% Obrázky v textu
\graphicspath{ {pictures/} } %% Cesta k obrázkům
\usepackage{tikz}
\usepackage{comment}

\usepackage{bbm}

\usepackage[
backend=biber,
style=alphabetic,
sorting=ynt
]{biblatex}
\title{When Mathematics Helps Physics: Calculation of the Integral of Kholodenko and Silagadze}
\author{Dominik Beck\\
\small \texttt{beckd@karlin.mff.cuni.cz}}
\date{31 March 2025}

\theoremstyle{definition}
\theoremstyle{plain}
\newtheorem*{theorem*}{Theorem}%theorem without numbering

\theoremstyle{remark}

\input{mycommands.sty}

\usepackage{titlesec}
\newcommand{\sectionfont}{\Large\bfseries}%Font size of sections headings
\titleformat{\section}
    {\titlerule     
     \vspace{0.3ex}%
     \sectionfont}
    {\thesection}{1em}
    {\sectionfont}[\titlerule\vspace{0.18em}]

\usepackage{listings,xcolor}

\newlength{\RoundedBoxWidth}
\newsavebox{\GrayRoundedBox}
\newenvironment{GrayBox}[1][\dimexpr\columnwidth-4.5ex]%
   {\setlength{\RoundedBoxWidth}{\dimexpr#1}
    \begin{lrbox}{\GrayRoundedBox}
       \begin{minipage}{\RoundedBoxWidth}}%
   {   \end{minipage}
    \end{lrbox}
    \begin{center}
    \begin{tikzpicture}
       \draw node[draw=black,fill=black!10,rounded corners,%
             inner sep=2ex,text width=\RoundedBoxWidth]%
             {\usebox{\GrayRoundedBox}};
    \end{tikzpicture}
    \end{center}}

\setlist[itemize]{left=1pt}

\begin{document}
\maketitle

\begin{abstract}
In this paper, we show a physics-free derivation of a Landau-Zener type integral introduced by Kholodenko and Silagadze.
%$$\textstyle I_n \!=\! \int_{-\infty}^\infty \!\int_{-\infty}^{x_1}\!\!\cdots\! \int_{-\infty}^{x_{2n-1}} \!\!\cos(x_1^2-x_2^2)\cos(x_3^2-x_4^2)\cdots \cos(x_{2n-1}^2-x_{2n}^2) \dd x_1 \dd x_2 \cdots \dd x_{2n} \!=\! \frac{2}{n!}\!\left(\frac{\pi}{4}\right)^{\!n}.$$
\end{abstract}

\bgroup
\itshape
\small
The study was supported by the Charles University, project GA UK No. 71224 and by Charles University Research Centre program No. UNCE/24/SCI/022.
\egroup

\tableofcontents

\section*{Introduction}
\addcontentsline{toc}{section}{Introduction}

In 2012, Kholodenko and Silagadze published an enigmatic paper \cite{kholodenko2012physics} in which they showed for $n = 1,2,3,\ldots $ that
\begin{equation}\label{Eq:In}
\begin{split}
    I_n \!=\! & \int_{-\infty}^\infty \!\int_{-\infty}^{x_1}\!\!\!\cdots\! \int_{-\infty}^{x_{2n-1}} \!\!\cos(x_1^2-x_2^2)\cos(x_3^2-x_4^2)\cdots \\
    & \cos(x_{2n-1}^2-x_{2n}^2) \ddd x_{2n} \cdots \dd x_2\dd x_1 \!=\! \tfrac{2}{n!}\!\left(\tfrac{\pi}{4}\right)^{\!n}\!\!.
\end{split}
\end{equation}
Their approach involved an intricate use of the concepts from quantum physics, namely the spinorization and the Hopf map, using which they were able to show that the integral is equivalent to the \emph{Landau-Zener formula}. The Landau-Zener formula expresses the exact probability amplitude of a certain quantum system with a dynamical potential and its standard derivation uses either the contour integration (Landau's approach) or asymptotics of a certain differential equation (Zener's approach), see Piquer i Méndez \cite{piquer2023adiabatic} for a clear exposition of both approaches. Alternatively, see Rojo \cite{rojo2010matrix} for a matrix exponential solution or Glasbrenner and Schleich \cite{glasbrenner2023landau} for a derivation using Markov property.

\vspace{0.5em}
We were aware of the result \eqref{Eq:In} as early as 2015, in which we solved an integral related to $I_4$ also introduced by Silagadze \cite{I4silagadze}. See our derivation here \cite{machinatoI4}. There is also an alternative derivation of the $I_4$ related integral due to de Reyna \cite{de2015cancellations} using clever asymptotics of a certain sine-type power series.

\vspace{0.5em}
However, back then, we lacked a proper mathematical background to tackle the general problem of $I_n$ using only simple and purely mathematical concepts as we did in the case of $I_4$. Now, exactly ten years later, we are returning to the problem with a new toolkit: Contour integration. Unlike the original Landau derivation, we will focus only on the integral without mentioning its connection to quantum physics. Despite its ubiquitousness and broad applications in quantum physics, simple derivations of Landau-Zener formula are rare, many of them use not so clear assumptions. The purpose of our paper is thus to present a physics-free derivation, mainly to make the result accessible to a broader mathematical community.

%\vspace{1em}
%Last but not least, we hope this paper and the derivation of $I_n$ herein will serve as a neat application of the well known techniques of the field of complex analysis.

\section{Preliminaries}
\subsection{Auxiliary functions}
In our paper, we define four auxiliary functions of a real parameter $t$ which turn out to be relevant,
\begin{equation}\label{Eq:UVPQ}
\begin{split}
    U(t) &= \int_0^\infty \frac{\sqrt{x}}{1+x^2} \exp\left(\frac{t}{4}\arctan x\right)\dd x, \\
    V(t) &= \int_0^\infty \frac{1/\sqrt{x}}{1+x^2} \exp\left(\frac{t}{4}\arctan x\right)\dd x,\\
    P(t) &= \int_0^1 \frac{\sqrt{y}-1}{1-y^2} \exp\left(\frac{it}{4}\argtanh y\right) \dd y,\\
    Q(t) &= \int_0^1 \frac{\frac{1}{\sqrt{y}}-1}{1-y^2} \exp\left(\frac{it}{4}\argtanh y\right) \dd y.
\end{split}
\end{equation}
Those functions can be expressed in an exact form using Appel and Generalized Hypergeometric functions. The exact expressions are, however, not necessary for our derivation of $I_n$. Note that those functions are not independent. By substitution $x \to 1/x$, we obtain an obvious relation between $U$ and $V$,
\begin{equation}
U(t) = e^{\pi t/8} \, V(-t).    
\end{equation}
Moreover, the complex-valued functions $P,Q$ can be expressed as linear combinations of real-valued-only functions $U$ and $V$. In order to find those relations, we need to study a certain contour integral.

\subsection{Contour integral}

\subsubsection{Definitions}
Consider the following contour integral with parameters $A,B \in \mathbb{C}$ and $t\in \mathbb{R}$,
\begin{equation}\label{Eq:contourint}
\oint_C \left[\frac{\sqrt{z}-\sqrt{i}}{1+z^2}A+ \frac{\frac{1}{\sqrt{z}}-\frac{1}{\sqrt{i}}}{1+z^2}B \right] e^{\frac{t}{4} \arctan z} \ddd z
\end{equation}
where we define $\arctan z = \frac{1}{2i}\ln(\frac{1+i z}{1-i z})$. Both $\sqrt{z}$ and $\ln z$ functions are assumed to have branch cuts at the negative real axis (the so called \emph{principal branch}), that is $\arg z \in (-\pi,\pi]$. Before proceeding further, let us clarify all the functions appearing in the contour integral. Readers familiar with complex functions may skip the following paragraphs. Writing $z = r e^{i\varphi}$ with $r>0$ and $ \varphi\in (-\pi,\pi]$, the complex square root function is defined by
\begin{equation}
    \sqrt{r e^{i\varphi}} = \sqrt{r} e^{i\varphi/2},
\end{equation}
so $\sqrt{i}$ in the formula above is given by
\begin{equation}
    \textstyle \sqrt{i} = e^{i\pi/4} = \frac{\sqrt{2}}{2} + \frac{\sqrt{2}}{2} i.
\end{equation}
Similarly, we have for the complex logarithm,
\begin{equation}
    \ln (re^{i\varphi}) = \ln r + i\varphi.
\end{equation}
Those definitions ensure both $\sqrt{z}$ and $\ln z$ are analytic on $\mathbb{C}\setminus (-\infty,0]$.

\subsubsection{Real axis}
Note that we recover the ordinary $\arctan$ function on reals. To see this, let $z = x,\, x \in \mathbb{R}$, then $\arg z \in (-\frac{\pi}{2},\frac{\pi}{2})$ and thus
\begin{equation}
    1 \pm ix = \sqrt{1+x^2} e^{\pm i \arctan x}.
\end{equation}
Substituting into the complex definition of $\arctan$,
\begin{equation}
    \textstyle \frac{1}{2i} \ln e^{2i\arctan x} = \frac{1}{2i}(2i\arctan x) = \arctan x
\end{equation}
since $2 \arctan x \in (-\pi,\pi)$.

\subsubsection{Branch cuts}
To locate the branch cuts of the complex function $\arctan z$, we need to find which points on the $z$-plane are mapped onto the negative real axis in $w = \frac{1+i z}{1 - iz}$ as $\arctan z = \frac{1}{2i}\ln w$. Solving
\begin{equation}
    \frac{1+iz}{1-iz} = -y, \quad y>0,
\end{equation}
we get the location of the branch cut in $z\in \mathbb{C}$ as
\begin{equation}
    z = i \frac{1+y}{1-y} \in (-\infty i,-i)\cup (i,i\infty),
\end{equation}
where $(z_1,z_2)$ denotes a complex interval from $z_1$ to $z_2$ (a line segment between any $z_1,z_2\in \mathbb{C}$).

\subsubsection{Imaginary axis}
Let $z=i y, y>0$, so $z = y e^{\pi i/2}$ and thus
\begin{equation}
    \sqrt{z} = \sqrt{i y} = \sqrt{y} e^{\pi i/4} =\sqrt{i} \sqrt{y}.
\end{equation}
The situation is different for the $\arctan$ function since it is no longer analytic on the whole imaginary axis. Let $z = i y, y\in(0,1)$ (analytic part), then
\begin{equation}
    \arctan z = \frac{1}{2i}\ln\left(\frac{1-y}{1+y}\right) = i\argtanh y,
\end{equation}
where $\argtanh y = \frac{1}{2}\ln\frac{1+y}{1-y}$ is the usual inverse hyperbolic tangent function. Right to the branch cut of $\arctan z$, we can write $z = i y + \varepsilon,y > 1$ and with $ \varepsilon >0$ small. By Taylor expansion,
\begin{equation}
\begin{split}
        \frac{1+iz}{1-iz}\bigg{|}_{iy+\varepsilon} \!\! &= \frac{1+iz}{1-iz}\bigg{|}_{iy}\!\!+\varepsilon\left(\frac{1+iz}{1-iz}\right)'\bigg{|}_{iy} \!\!\!+ O(\varepsilon^2)\\
        &= -\frac{y-1}{y+1} + \frac{2i\varepsilon}{(y-1)^2} = \frac{y-1}{y+1} e^{\pi i},
\end{split},
\end{equation}
which gives
\begin{equation}
    \textstyle\arctan z = \frac{1}{2i}\left(\ln\frac{y-1}{y+1}+\pi i\right)=\frac{\pi}{2}+i\argtanh\frac{1}{y}.
\end{equation}

\subsubsection{Asymptotics}
Let $|z| \to \infty$, expanding $\frac{1+iz}{1-iz}$, we get
\begin{equation}
\textstyle \frac{1+iz}{1-iz} = -1 + \frac{2i}{z} + O(\frac{1}{z^2}) = e^{\pi i}\left(1 - \frac{2i}{z}\right) + O(\frac{1}{z^2}).
\end{equation}
If $\Re z \to \infty$ only (that is $z \to \infty e^{i\varphi}$ with $\varphi \in (-\frac{\pi}{2},\frac{\pi}{2})$), we get $\arg (1-2i/z)\in (-\pi,0)$ and thus
\begin{equation}
\textstyle\ln\left(\frac{1+iz}{1-iz}\right) \!=\! \ln(e^{\pi i})+\ln\left(1\!-\!\frac{2i}{z}\right)=\pi i - \frac{2i}{z} + O(\frac{1}{z^2}),
\end{equation}
from which, dividing by $2i$,
\begin{equation}
    \arctan z = \frac{\pi}{2} - \frac{1}{z} + O\left(\frac{1}{z^2}\right), \quad \Re z \to \infty.
\end{equation}

\subsubsection{Pole expansion}
Let us examine the behaviour of $\arctan z$ near its singularity $z=i$ in the first quadrant in $\mathbb{C}$. Let $z = i - i \varepsilon e^{i \varphi}, \varepsilon>0$ and $\varphi \in (0,\pi)$. Then
\begin{equation}
    \frac{1+iz}{1-iz} = \frac{1-1+\varepsilon e^{i\varphi}}{1+1-\varepsilon e^{i\varphi}} = \frac{\varepsilon}{2} e^{i \varphi} + O(\varepsilon^2),
\end{equation}
from which
\begin{equation}
    \arctan z = \frac{1}{2i} \ln \frac{1+iz}{1-iz} = \frac{1}{2i}\ln\frac{\varepsilon}{2} + \frac{\varphi}{2} + O(\varepsilon^2). 
\end{equation}
Next, the square brackets part of $f(z)$ has a finite limit as $z \to i$. By L'Hospitals rule,
\begin{equation}
    \lim_{z\to i} \frac{\sqrt{z}-\sqrt{i}}{1+z^2} = \frac{\frac{1}{2\sqrt{z}}}{2z}\bigg{|}_{z=i}=\frac{1}{4i\sqrt{i}},
\end{equation}
similarly for $(\frac{1}{\sqrt{z}}-\frac{1}{\sqrt{i}})/(1+z^2)$ as $z \to i$.

\subsection{Parametrisation}
\begin{figure}[h]
    \centering
    \includegraphics[width=0.6\linewidth]{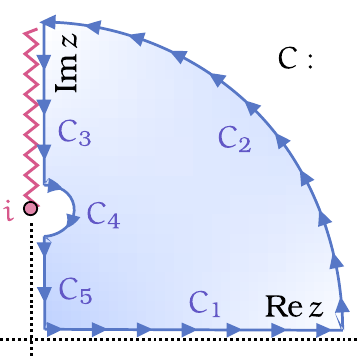}
    \caption{Contour $C$}
    \label{Fig:Contour}
\end{figure}
We are now ready to examine the integral \eqref{Eq:contourint}. It is convenient to denote
\begin{equation}
    f(z) = \left[\frac{\sqrt{z}-\sqrt{i}}{1+z^2}A+ \frac{\frac{1}{\sqrt{z}}-\frac{1}{\sqrt{i}}}{1+z^2}B \right] e^{\frac{t}{4} \arctan z}.
\end{equation}
Let us consider the counter-clockwise contour $C$ consisted of integration curves $C_k$, $k=1,\ldots,6$,
\begin{equation}
    C = C_1 \cup C_2 \cup C_3 \cup C_4 \cup C_5
\end{equation}
as shown in Figure \ref{Fig:Contour}. 

\vspace{1em}
Inside $C$, the function $f$ is analytic, so by the Cauchy integral formula,
\begin{equation}
    \oint_C f(z) \dd z = \sum_{k=1}^6 \int_{C_k} f(z)\dd z =0.
\end{equation}
Let us now parametrize the integrals on the individual integration curves, we denote $J_k = \int_{C_k} f(z) \dd z$ and $\circleddash C_k$ as negatively oriented $C_k$.
\begin{itemize}
    \item $C_1: z=x,\quad x\in (0,R),\quad R \to \infty,\quad \dd z = \dd x$,
    \begin{equation}
        J_1 \!=\! \int_0^R \!\left[\!\frac{\sqrt{x}\!-\!\sqrt{i}}{1+x^2}A\!+\! \frac{\frac{1}{\sqrt{x}}\!-\!\frac{1}{\sqrt{i}}}{1+x^2}B \right]\! e^{\frac{t}{4} \arctan x}\ddd x.
    \end{equation}
    As $R\to \infty$, we may express the limit using the auxiliary functions defined earlier,
    \begin{equation}
        \!\! J_1 \!\!\to\! A U(t) + B V(t) -\! \sqrt{i}(A-iB) \!\!\int_0^\infty \!\!\frac{e^{\frac{t}{4}\!\arctan x}}{1+x^2}\dd x.
    \end{equation}
    The integral on the right is trivial, we get
    \begin{equation}
        J_1 \!\to\! A U(t) + B V(t) - \frac{4 \sqrt{i}}{t} (A-iB)\left(e^{\frac{\pi t}{8}}\!-\!1\right).
    \end{equation}
    \item $C_2: z=R e^{i\varphi}, R\to \infty, \varphi \in (0,\frac{\pi}{2}), \dd z = iR e^{i\varphi} \dd \varphi,$
    \begin{equation}
        |J_2| \leq R \int_0^{\frac{\pi}{2}} \!\left[\!\frac{\sqrt{R} +\! 1}{R^2\!-\!1}|A|\!+\! \frac{\frac{1}{\sqrt{R}} \!+\! 1}{R^2\!-\!1}|B| \right] \! e^{\frac{\pi |t|}{8}}\dd \varphi,
    \end{equation}
    from which $J_2 \to 0$ as $R \to \infty$.
    \item $\circleddash C_3: z=iy + \varepsilon,\quad y \in (1+\varepsilon,R),\quad \dd z = i \dd y,$
    \begin{equation}
        \begin{split}
        J_3 = &-i\sqrt{i}\int_{1-\varepsilon}^R \exp\left(\frac{\pi t}{8} +\frac{i t}{4}\argtanh y\right) \\
        & \times\bigg{[}\frac{\sqrt{y}-1}{1-y^2}A-i \frac{\frac{1}{\sqrt{y}}-1}{1-y^2} B \bigg{]} \dd y.
        \end{split}
    \end{equation}
    Substituting $y\to 1/y$ and letting $R\to \infty$ and $\varepsilon \to 0^+$, we get, using auxiliary functions
    \begin{equation}
        J_3 \to i\sqrt{i} e^{\frac{\pi t}{8}}(A Q(t)-i B P(t))
    \end{equation}
    \item $\circleddash C_4: z\!=\!i-i\varepsilon e^{i\varphi}, \varphi \in (0,\pi), \varepsilon\!\to\! 0^+\!, \dd z \!= \!\varepsilon e^{i\varphi}\dd \varphi,$
    \begin{equation}
        |J_4| \leq \varepsilon \int_0^\pi \left[\frac{|A|}{4}+\frac{|B|}{4}+O(\varepsilon)\right] e^{\frac{t}{8i}\ln\frac{\varepsilon}{2}+\frac{\varphi t}{8}} \dd \varphi.
    \end{equation}
    Since $t\in \mathbb{R}$, $J_4$ vanishes as $\varepsilon \to 0^+$.
    \item $\circleddash C_5: z=iy, y \in (1,1-\varepsilon), \varepsilon\to 0^+, \dd z = i \dd y,$
    \begin{equation}
        \textstyle \!\!\!\! J_5 \!=\! -i\sqrt{i}\!\int_0^{\!1\!-\!\varepsilon} \!\!\left[\!\frac{\sqrt{y}-1}{1-y^2}A\!-\!i \frac{\frac{1}{\sqrt{y}}-1}{1-y^2} B \right]\! e^{\!\frac{i t}{4}\!\argtanh y}\dd y.\!\!
    \end{equation}
    As $\varepsilon\to 0^+$, we get, using auxiliary functions,
    \begin{equation}
        J_5 \to -i\sqrt{i} \left(A P(t)-i B Q(t)\right).
    \end{equation}
\end{itemize}

\subsection{Relations connecting the auxiliary functions}
By the Cauchy integral formula, $J_1 + \cdots + J_5 = 0$. Taking the limit as $R \to \infty$ and $\varepsilon \to 0^+$, we get for any $A,B \in \mathbb{C}$,
\begin{equation}
\begin{split}
    & \! A U(t) + B V(t) = \frac{4\sqrt{i}}{t}(A-iB)(e^{\frac{\pi t}{8}}-1)\\
    & \!+i\sqrt{i}\left(A P(t)\!-\!i B Q(t)\!-\!e^{\frac{\pi t}{8}}AQ(t)\!+\!i e^{\frac{\pi t}{8}}B P(t)\right).
\end{split}
\end{equation}
Comparing the terms with the coefficient $A$ and $B$, respectively, we get
\begin{equation}
\begin{split}
        U(t) & = \frac{4\sqrt{i}}{t}(e^{\frac{\pi t}{8}}-1)+i\sqrt{i}\left(P(t)-e^{\frac{\pi t}{8}}Q(t)\right),\\
        V(t) & =\! -\frac{4i\sqrt{i}}{t}(e^{\frac{\pi t}{8}}\!-\!1)+\sqrt{i}\left(Q(t)\!-\!e^{\frac{\pi t}{8}}P(t)\right).
\end{split}
\end{equation}
Solving for $P$ and $Q$, we get, finally,
\begin{equation}\label{Eq:PQreal}
    \begin{split}
        P(t) & = -\frac{4 i}{t}+\sqrt{i}\,\,\frac{U(t)+i e^{\frac{\pi  t}{8}} V(t)}{e^{\frac{\pi  t}{4}}-1},\\
        Q(t) &= -\frac{4 i}{t}+\sqrt{i}\,\,\frac{e^{\frac{\pi  t}{8}} U(t)+i V(t)}{e^{\frac{\pi  t}{4}}-1}.
    \end{split}
\end{equation}

\section{Derivation}
\subsection{Integral equation formulation}
The relation \eqref{Eq:In} can be written recursivelly by introducing new functions $\tau_n$ as
\begin{equation}\label{Eq:taun}
    \tau_n(x) = \int_{-\infty}^x \int_{-\infty}^{x_1} \cos(x_1^2-x_2^2)\, \tau_{n-1}(x_2) \ddd x_2 \dd x_1
\end{equation}
with initial condition $\tau_0(x) = 1$. We have
\begin{equation}
    I_n = \tau_n(\infty), \quad n=1,2,3,\ldots.
\end{equation}
Let us introduce a new parameter $t\in \mathbb{R}$ using which we define a generating function $T$ as
\begin{equation}
    T(x,t) = \sum_{n=0}^\infty \tau_n(x) t^n.
\end{equation}
Equation \eqref{Eq:taun} then gets transformed into
\begin{equation}\label{Eq:Txtintegraleq}
    \!T(x,t) \!= \!1  + t \int_{-\infty}^x \int_{-\infty}^{x_1} \!\!\cos(x_1^2-x_2^2)\, T(x_2,t) \dd x_2 \dd x_1,
\end{equation}
which is an integral equation for $T(x,t)$ subject to the boundary condition
\begin{equation}
    T(-\infty,t) = 1.
\end{equation}
The solution of \eqref{Eq:Txtintegraleq} at $x \to \infty$ completely solves our problem of finding $I_n$ since by expanding $T(\infty,t)$ in $t$ when $t \to 0$, we get
\begin{equation}\label{Eq:Tseries}
    T(\infty,t) = \sum_{n=0} \tau_n(\infty) t^n = 1 + \sum_{n=1}^\infty I_n t^n.
\end{equation}
Note that, rescaling by $\sqrt{s}$ with $s>0$ and letting $s \to \infty$, we get the \emph{stretching relations} of $T(x,t)$,
\begin{equation}\label{Eq:StrechRel}
\begin{split}
        \lim_{s \to \infty} T(\tfrac{x}{\sqrt{s}},t) & = T(0,t),\\
        \lim_{s \to 0^+} T(\tfrac{x}{\sqrt{s}},t) & = T(-\infty,t)\mathbbm{1}_{x<0} + T(\infty,t)\mathbbm{1}_{x>0},
\end{split}
\end{equation}
we will see their importance later.

\subsection{Differential equation formulation}
It is convenient to denote
\begin{equation}\label{Eq:alphabeta}
    \begin{split}
        \alpha(t) & = \int_{-\infty}^0 \cos (x^2) T(x,t)\ddd t,\\ 
        \beta(t) & = \int_{-\infty}^0 \sin (x^2) T(x,t)\ddd t.
    \end{split}
\end{equation}
Differentiating the integral equation \eqref{Eq:Txtintegraleq} with respect to $x$, we get,
\begin{equation}\label{Eq:dTdx}
\frac{\partial T}{\partial x} = t \int_{-\infty}^x  \cos(x^2-y^2) T(y,x) \ddd y.
\end{equation}
At $x=0$, we get $\frac{\partial T}{\partial x}(0,t) = t\alpha(t)$. Second differentiation with respect of $x$ gives
\begin{equation}
    \frac{\partial^2 T}{\partial x^2} = t T(x,t) -2xt \int_{-\infty}^x \sin(x^2-y^2)T(y,t) \ddd y.
\end{equation}
At $x=0$, we get $\frac{\partial^2 T}{\partial x^2}(0,t) = t T(0,t)$. Yet another differentiation yields
\begin{equation}
    \frac{\partial^3 T}{\partial x^3}\!=\!t\frac{\partial T}{\partial x} -4x^2 \frac{\partial T}{\partial x} \!-\! 2t\!\int_{-\infty}^x \!\!\sin(x^2\!-\!y^2)T(y,t)\dd y,
\end{equation}
which gives $\frac{\partial^3 T}{\partial x^3}(0,t) = t^2 \alpha(t)-2t\beta(t)$ at $x=0$. One final differentiation results in the following differential equation
\begin{equation}
    \frac{\partial^4 T} {\partial x^4} =(t- 4x^2)\frac{\partial^2 T}{\partial x^2} - 12 x \frac{\partial T}{\partial x}.
\end{equation}

This is an ordinary differential equation (ODE) with parameter $t$. Its solution at $x \to \infty$ completely solves our problem of finding $I_n$ (expanding $T(\infty,t)$ in $t$). For completeness, the following is the complete list of \emph{boundary conditions} on $T(x,t)$,
\begin{align*}
        T(-\infty,t) & = 1, && \frac{\partial^2 T}{\partial x^2}(0,t) = t T(0,t),\\
        \frac{\partial T}{\partial x}(0,t) &= t \alpha(t),
        && \frac{\partial^3 T}{\partial x^3}(0,t) = t^2 \alpha(t)-2t\beta(t).
\end{align*}

\begin{figure}[h]
    \centering
    \includegraphics[width=1.0\linewidth]{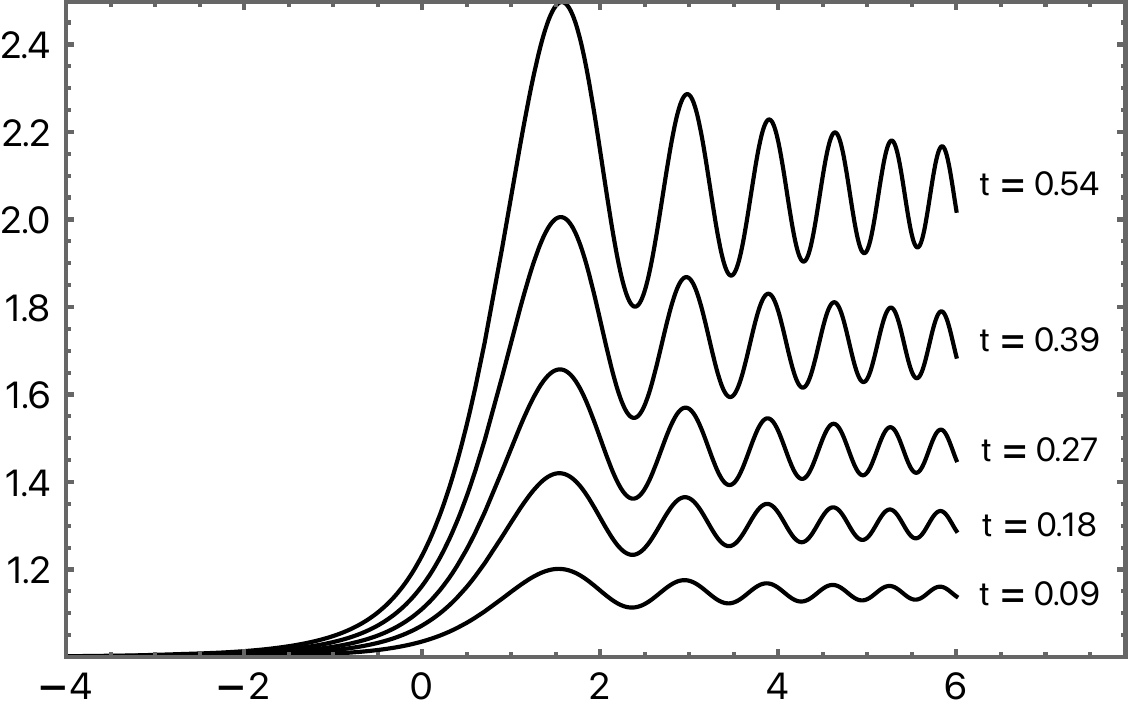}
    \caption{Numerical solution of ODE for $T(x,t)$}
    \label{Fig:Sample}
\end{figure}

Using Mathematica $\mathtt{NDSolve}$ command, we can see numerically how the solution looks like for various parameters $t$ (Figure \ref{Fig:Sample}).

\subsection{Even integral transform}
In what follows, we establish a relation between $T(\infty,t)$ and $T(0,t)$. Let
\begin{equation}
    e(s,t) = \int_{-\infty}^\infty T(x,t)\, e^{-sx^2}\ddd x.
\end{equation}
Since $T(x,t)$ is continuous, bounded and having a finite limit $T(\infty,t)$, the integral on the right converges for any $s\in \mathbb{C}$ with $\Re s > 0$. Assuming $s$ is a positive real number, we get by substitution $x= y/\sqrt{s}$,
\begin{equation}
    e(s,t) = \frac{1}{\sqrt{s}}\int_{-\infty}^\infty T(\tfrac{y}{\sqrt{s}},t) e^{-y^2} \ddd y
\end{equation}
It is thus convenient to define $E(s,t) = e(s,t)\sqrt{s}$. Imposing the stretching relations \eqref{Eq:StrechRel},
\begin{align}
        \!\!E(\infty,t) & \!=\! \int_{-\infty}^\infty T(0,t) e^{-y^2}\ddd y = \sqrt{\pi}\,T(0,t),\!\\ 
\begin{split}\label{Eq:E0plus}
        \!\!E(0^+\!\!,t) & \!=\! \int_{-\infty}^\infty \!\!\!\!\! (T(-\infty,t)\mathbbm{1}_{y<0} \! +\! T(\infty,t) \mathbbm{1}_{y>0} )e^{-y^2}\!\dd y \\ & \!= \frac{\sqrt{\pi}}{2}(T(\infty,t)+1).
\end{split}
\end{align}
Differentiating $e(s,t)$ with respect to $s$ and integrating by parts with respect to $x$ yields
\begin{equation}
    -e(s,t)-2s\frac{\partial e}{\partial s} = \int_{-\infty}^\infty x \frac{\partial T}{\partial x} e^{-sx^2}\ddd x
\end{equation}
or in terms of $E(s,t)$,
\begin{equation}
    -2\sqrt{s}\frac{\partial E}{\partial s} = \int_{-\infty}^\infty x \frac{\partial T}{\partial x} e^{-sx^2}\ddd x.
\end{equation}
Writing out the right hand side using Equation \eqref{Eq:dTdx} and changing the order of integration,
\begin{equation}
    \!-2\sqrt{s} \frac{\partial E}{\partial s} \!=\! t \int_{-\infty}^\infty \!\!\!\! T(y,t) \int_y^\infty \!\!\cos(x^2\!-\!y^2) x e^{-sx^2} \dd x \dd y.
\end{equation}
The inner integral is trivial and is equal to $\frac{s}{2(1+s^2)}e^{-sy^2}$. Hence, we get a differential equation for $E(s,t)$
\begin{equation}
    \frac{\partial E}{\partial s} = -\frac{t E(s,t)}{4(1+s^2)}.
\end{equation}
The solution of this ODE with the appropriate boundary condition $E(\infty,t)=\sqrt{\pi}\,T(0,t)$ must be
\begin{equation}
    E(s,t) = \sqrt{\pi}\, T(0,t) \exp\left(\frac{t}{4}\arccot s\right).
\end{equation}    
The remaining condition on $E(0^+,t)$ (Equation \eqref{Eq:E0plus}) therefore yields as $s \to 0^+$,
\begin{GrayBox}
\begin{equation}\label{Eq:Tinfty}
    T(\infty,t) = 2T(0,t)\exp\left(\tfrac{\pi t}{8}\right) - 1.
\end{equation}    
\end{GrayBox}

\subsection{Truncated integral transform}
In this section, we establish a relation between $T(-\infty,t)=1$ and $T(0,t)$. Let
\begin{equation}
    d(s,t) = \int_{-\infty}^\infty T(x,t)\, e^{-sx^2} \mathbbm{1}_{x<0}\ddd x.
\end{equation}
That is, the integration domain is $(-\infty,0)$. Again, the integral converges for any $s\in \mathbb{C}$ with $\Re s > 0$. Assuming $s$ is a positive real number, we get by substitution $x= y/\sqrt{s}$,
\begin{equation}
    d(s,t) = \frac{1}{\sqrt{s}}\int_{-\infty}^0 T(\tfrac{y}{\sqrt{s}},t) e^{-y^2} \ddd y
\end{equation}
It is thus convenient to define analogously $D(s,t) = d(s,t)\sqrt{s}$. Imposing the stretching relations \eqref{Eq:StrechRel},
\begin{align}
        D(\infty,t) & = \int_{-\infty}^0 T(0,t) e^{-y^2}\ddd y = \frac{\sqrt{\pi}}{2}\, T(0,t),\label{Eq:Dinfty}\\ 
        D(0^+,t) & = \int_{-\infty}^0  T(-\infty,t) e^{-y^2}\dd y = \frac{\sqrt{\pi}}{2}.
\end{align}
Differentiating $d(s,t)$ with respect to $s$ and integrating by parts with respect to $x$ yields
\begin{equation}
    -d(s,t)-2s\frac{\partial d}{\partial s} = \int_{-\infty}^\infty x \frac{\partial T}{\partial x} e^{-sx^2}\mathbbm{1}_{x<0}\ddd x
\end{equation}
or in terms of $D(s,t)$,
\begin{equation}
    -2\sqrt{s}\frac{\partial D}{\partial s} = \int_{-\infty}^\infty x \frac{\partial T}{\partial x} e^{-sx^2} \mathbbm{1}_{x<0} \ddd x.
\end{equation}
Writing out the right hand side using Equation \eqref{Eq:dTdx} and changing the order of integration,
\begin{equation}
    \!-2\sqrt{s} \frac{\partial D}{\partial s} \!=\! t \int_{-\infty}^0 \!\!\!\! T(y,t) \int_y^0 \!\!\cos(x^2\!-\!y^2) x e^{-sx^2} \dd x \dd y.
\end{equation}
The inner integral is again trivial and is equal to
\begin{equation}
    \frac{se^{-sy^2}}{2(1+s^2)} - \frac{s \cos (y^2)+\sin(y^2)}{2(1+s^2)}.
\end{equation}
Hence, using the definition of $\alpha(t)$ and $\beta(t)$ (Equation \eqref{Eq:alphabeta}), we get a differential equation for $D(s,t)$
\begin{equation}
    \frac{\partial D}{\partial s} = -\frac{t D(s,t)}{4(1+s^2)} + \frac{t}{4\sqrt{s}}\left(\frac{s\alpha(t)}{1+s^2}+\frac{\beta(t)}{1+s^2}\right).
\end{equation}
The solution of this ODE with the boundary condition $D(0^+,t)=\sqrt{\pi}/2$ is
\begin{equation}\label{Eq:solDst}
\begin{split}
    D(s,t) & = \frac{\sqrt{\pi}}{2}\, e^{-\frac{s}{4}\arctan s}\bigg{[}1+\frac{t}{2\sqrt{\pi}}\int_0^s \frac{1}{\sqrt{r}}\\
    & \times\left(\frac{r \alpha(t)}{1+r^2}+\frac{\beta(t)}{1+r^2}\right) e^{\frac{t}{4}\arctan r}\dd r\bigg{]}.
\end{split}
\end{equation}    
The remaining condition on $D(\infty,t)$ (Equation \eqref{Eq:Dinfty}) therefore yields, plugging $s \to \infty$,
\begin{GrayBox}
\vspace{-1em}
\begin{equation}\label{Eq:T0}
    \textstyle T(0,t) e^{\frac{\pi t}{8}} = 1 + \frac{t}{2\sqrt{\pi}}\left(\alpha(t) U(t)+\beta(t) V(t)\right),
\end{equation}    
\end{GrayBox}
where $U$ and $V$ are defined by Equation \eqref{Eq:UVPQ}.

\subsection{Imaginary boundary conditions}
The key to our problem is the realization that the solution we found for $D(s,t)$ is valid for all $s \in \mathbb{C}$ for which $\Re s >0$. The functions $\alpha(t)$ and $\beta(t)$ are not chosen arbitrarily since they are defined using $T(x,t)$. We may express them (Eqution \eqref{Eq:alphabeta}) also in terms of $d(s,t)$ as
\begin{equation}
\begin{split}
    \alpha(t) \mp i \beta(t) & =\int_{-\infty}^0 T(x,t)(\cos x^2\mp i\sin x^2)\dd x \\
    & = \int_{-\infty}^0 T(x,y) e^{\mp i x^2} \dd x = d(\pm i,t).
\end{split}
\end{equation}
The value $d(\pm i,t)$ must be viewed as a limit of $d(s,t)$ for $s \to \pm i$ (since $\Re s >0$). This gives us the final boundary condition. Considering only the case $s \to i$, we get in terms of $D(s,t) = d(s,t)\sqrt{s}$,
\begin{equation}
    D(i,t) = \sqrt{i}\, (\alpha(t)-i \beta(t)).
\end{equation}
On the other hand, $D(i,t)$ can be calculated directly from the solution we found for $D(s,t)$. Let $s=i\sigma, \sigma \in (0,1)$, then by Equation \eqref{Eq:solDst},
\begin{equation}
\begin{split}
    &D(i \sigma,t) =\frac{\sqrt{\pi}}{2}e^{-\frac{it}{4}\argtanh\sigma}\bigg{[}1+\frac{i\sqrt{i} \,t}{2\sqrt{\pi}}\times \\
    & \int_0^\sigma\!\left(\!\frac{\sqrt{y}}{1\!-\!y^2}\alpha(t)-\frac{ \frac{1}{\sqrt{y}}}{1\!-\!y^2}i\beta(t)\!\right) e^{-\frac{it}{4}\argtanh\sigma} \ddd y\bigg{]}.
\end{split}
\end{equation}
Writing
\begin{equation}
\textstyle\frac{\sqrt{y}}{1-y^2} = \frac{\sqrt{y}-1}{1-y^2} + \frac{1}{1-y^2},\quad \frac{\frac{1}{\sqrt{y}}}{1-y^2} = \frac{\frac{1}{\sqrt{y}}-1}{1-y^2} + \frac{1}{1-y^2}
\end{equation}
and integrating out the $1/(1-y^2)$ term, we get
\begin{equation}
\begin{split}
    & D(i\sigma,t) = \sqrt{i}(\alpha(t)-i\beta(t)) + \frac{\sqrt{\pi}}{2}e^{-\frac{it}{4}\argtanh\sigma}\\
    &\times\bigg{[}1-\frac{2\sqrt{i}}{\sqrt{\pi}}(\alpha(t)-i\beta(t))+\frac{i\sqrt{i} \,t}{2\sqrt{\pi}}\times \\
    & \int_0^\sigma\!\!\left(\!\frac{\sqrt{y}\!-\!1}{1\!-\!y^2}\alpha(t)-\frac{ \frac{1}{\sqrt{y}}\!-\!1}{1\!-\!y^2}i\beta(t)\!\right)\! e^{-\frac{it}{4}\argtanh\sigma} \dd y\bigg{]}
\end{split}
\end{equation}

Note that the integral is now finite as $\sigma \to 1$. Crucially, as $\sigma$ approaches $1$, since the first part already yields the correct limit $\sqrt{i}(\alpha(t)-i\beta(t))$, the term in the square bracket must vanish, otherwise the limit as a whole would not exist. Hence, we get the following condition
\begin{equation}
\textstyle 1\!=\!\frac{2\sqrt{i}}{\sqrt{\pi}}(\alpha(t)-i\beta(t))-\frac{i\sqrt{i} \,t}{2\sqrt{\pi}}(\alpha(t)P(t)-i\beta(t)Q(t)),
\end{equation}
where $P$ and $Q$ are defined by Equation \eqref{Eq:UVPQ}. Substituting for $P$ and $Q$ from Equation \eqref{Eq:PQreal}, we get
\begin{equation}
\begin{split}
\frac{2 \sqrt{\pi }}{t}\left(e^{\frac{\pi t}{4}}-1\right) & =\alpha(t)U(t) + \beta(t)V(t) \\
& +i e^{\frac{\pi  t}{8}} (\alpha(t)V(t) - \beta(t)U(t) ).
\end{split}
\end{equation}
Comparing the real and imaginary part, we obtain
\begin{equation}
\begin{split}
\alpha(t)U(t) + \beta(t)V(t) & = \frac{2 \sqrt{\pi }}{t}\left(e^{\frac{\pi t}{4}}-1\right), \\
\alpha(t)V(t) - \beta(t)U(t) & = 0,
\end{split}
\end{equation}
which is a linear system of equations for $\alpha(t)$ and $\beta(t)$ with solution given as
\begin{comment}
\begin{equation}
    \!\!\!\!\!\alpha(t) \!=\! \frac{2 \sqrt{\pi } (e^{\frac{\pi  t}{4}}\!-\!1) U(t)}{t \left(U(t)^2+V(t)^2\right)},
    \beta(t) \!=\! \frac{2 \sqrt{\pi }
   (e^{\frac{\pi  t}{4}}\!-\!1) V(t)}{t \left(U(t)^2+V(t)^2\right)}.\!\!\!\!\!
\end{equation}
\end{comment}
\begin{equation}
\begin{split}
    \alpha(t) &= \frac{2 \sqrt{\pi } (e^{\frac{\pi  t}{4}}-1) U(t)}{t \left(U(t)^2+V(t)^2\right)},\\
    \beta(t) & = \frac{2 \sqrt{\pi }
   (e^{\frac{\pi  t}{4}}-1) V(t)}{t \left(U(t)^2+V(t)^2\right)}.
\end{split}
\end{equation}

\subsection{Final matching}
Substituting $\alpha(t)$ and $\beta(t)$ into Equation \eqref{Eq:T0},
\begin{equation}
    T(0,t) = \exp(\tfrac{\pi t}{8}),
\end{equation}
from which, by Equation \eqref{Eq:Tinfty},
\begin{equation}
    T(\infty,t) = 2 \exp(\tfrac{\pi t}{4}) - 1.
\end{equation}
The series expansion of $T(\infty,t)$ (Equation \eqref{Eq:Tseries}) immediately yields
\begin{equation}
    I_n = \tfrac{2}{n!}\left(\tfrac{\pi}{4}\right)^n.
\end{equation}

\section{Related result}
A similar analysis shows that also
\begin{equation}\label{Eq:In}
\begin{split}
    & K_n \!=\! \int_{-\infty}^\infty \!\int_{x_1}^\infty\!\!\!\cdots\! \int_{x_{2n-1}}^\infty \!\!\cos(x_1^2-x_2^2)\cos(x_3^2-x_4^2)\cdots \\
    & \quad\cos(x_{2n-1}^2-x_{2n}^2) \ddd x_{2n} \cdots \dd x_2\dd x_1 \!=\! \tfrac{2}{n!}\!\left(\tfrac{\pi}{4}\right)^{\!n}\!\!.
\end{split}
\end{equation}

\section*{Acknowledgements}
\addcontentsline{toc}{section}{Acknowledgements}
First of all I would like to thank professor Milo\v{s} Zahradn\'{i}k for sparking my interest in peculiar complex unsolved problems in probability theory during his insightful lectures held in Prague. I am also thankful to Lucie Abigail Kopelentov\'{a} for her fruitful insights and especially for her critical first careful reading and language corrections of my work.

\printbibliography[
heading=bibintoc]

\end{document}